\documentclass[10pt,leqno]{article}
\AtBeginDvi{} 
\usepackage{amssymb}
\usepackage{amscd}
\usepackage{amsmath}
\usepackage{amsthm}
\usepackage{mathrsfs}
\usepackage{graphics}
\usepackage{graphicx}
\usepackage{mediabb}
\def\overset#1#2{{\mathrel{\mathop {{#2}_{}}\limits^{#1}}}}
\def\underset#1#2{{\mathrel{\mathop {{}_{} {#2}}\limits_{{#1}_{}}}}}
\def\upplim_#1{\underset{#1}{\overline\lim}\;}
\def\lowlim_#1{\underset{#1}{\underline\lim}\;}
\def\ang#1{{\langle}#1{\rangle}}
\newtheorem{defn}[equation]{\indent{\it Definition}\rm }

\newtheorem{lem}[equation]{Lemma}

\newtheorem{rmk}[equation]{\indent \rm {\it Remark}}

\newtheorem{thm}[equation]{Theorem}

\newcommand{\C}{{\mathbf{C}}}
\newcommand{\sC}{\mathscr{C}}

\newcommand{\del}{{\partial}}
\newcommand{\delbar}{\bar{\partial}}

\newcommand{\ex}{{}^\exists}
\newcommand{\fa}{^\forall}
\newcommand{\sF}{\mathscr{F}}
\newcommand{\iso}{\cong}

\newcommand{\sI}{\mathscr{I}}
\newcommand{\im}{{\Im\,}}

\newcommand{\sM}{\mathscr{M}}

\renewcommand{\O}{{\mathcal{O}}}

\newcommand{\rP}{{\mathrm{P}}}
\newcommand{\PD}{{\mathrm{P}\Delta}}

\newcommand{\R}{{\mathbf{R}}}
\newcommand{\sR}{\mathscr{R}}
\newcommand{\re}{{\Re\,}}

\newcommand{\Sig}{\mathit{\Sigma}}

\newcommand{\sU}{\mathscr{U}}
\newcommand{\Z}{\mathbf{Z}}
\newenvironment{pf}{\par{\it Proof.\,}}{\qed\vskip+2pt}
\numberwithin{equation}{section}

\title{
A Weak Coherence Theorem and Remarks to the Oka Theory
}
\date{\empty}
\author{
By Junjiro Noguchi\thanks{Research supported in part by Grant-in-Aid
 for Scientific Research (C) 15K04917. \hfill\break
\hbox{\quad} AMC2010: 32A99; 32E30 \hfill\break
\hbox{\quad} Key words: coherence, Oka, Levi problem,  Hartogs' inverse problem,
 several complex variables
 \hfill\break
\hbox{\quad} Affiliation and address: {Graduate School of Mathematical Sciences},
{University of Tokyo (Emeritus)};
{Komaba, Meguro-ku, Tokyo 153-8914}, {Japan}\hfill\break
\hbox{\quad} e-mail: noguchi@ms.u-tokyo.ac.jp
}}
\begin{document}
\AtBeginDvi{} 
%
\setlength{\baselineskip}{16pt}
\parskip+2.5pt
\maketitle
\begin{abstract}
The proofs of K. Oka's Coherence Theorems are based on Weierstrass'
 Preparation (division) Theorem.
Here we formulate and prove a Weak Coherence Theorem
without using Weierstrass' Preparation Theorem,
but only with \textit{power series expansions}:
The proof  is almost of linear algebra.
Nevertheless, this simple {\em Weak Coherence Theorem} suffices to
give other proofs of 
the Approximation,  Cousin I/II, 
and Levi's (Hartogs' Inverse) Problems even in  simpler ways than those known,
 as far as the domains are
 non-singular; they constitute the main basic part of the theory
of several complex variables. 

The new  approach enables us to complete
 the proofs of those problems
in quite an elementary way
 without Weierstrass' Preparation Theorem or
 the cohomology theory of Cartan--Serre, nor $L^2$-$\delbar$ method 
of H\"ormander.

We will also recall some new historical facts that Levi's (Hartogs' Inverse) Problem
of general dimension $n \geq 2$ was, in fact,  solved by K. Oka in 1943
(unpublished) and by S. Hitotsumatsu in 1949 (published in Japanese),
whereas it has been usually recognized as proved by K. Oka 1953,
by H.J. Bremermann and by F. Norguet 1954, independently.
\end{abstract}

\section{Introduction and a weak coherence}

K. Oka \cite{ok7}, \cite{ok8} proved three fundamental coherence
theorems for
\begin{enumerate}\setlength{\itemsep}{-3pt}
\item
 the sheaf $\O:=\O_{\C^n}$ of germs of holomorphic functions on $\C^n$,
\item
 the ideal sheaf $\sI\ang{A}$ of an analytic subset $A$
of an open subset of $\C^n$,
\item
 the normalization of the structure sheaf of a complex space,
\end{enumerate}
where for the second, H. Cartan \cite{ca50} gave his own proof based on
 Oka \cite{ok7} (cf.\ \cite{nog16} Chap.~9).
We simply call $\sI\ang{A}$ a (resp.\ non-singular)
 {\em geometric ideal sheaf}
of a (resp.\ non-singular) analytic subset $A$ (cf., e.g., \cite{nog16}
Chap.~6).
Oka's Coherence has played a fundamental role in modern Mathematics,
so that it led to the notion of ringed spaces due to H. Cartan,
and developed by J.-P. Serre, R. Remmert,  H. Grauert and
A. Grothendieck (cf.\ \cite{dieu} p.~162).
The proofs of the above coherence theorems rely on Weierstrass'
Preparation (division) Theorem. 

The aim of this paper is to formulate 
 a {\em Weak Coherence Theorem}
(Theorem~\ref{weakc} below), which we prove
 \textit{not with Weierstrass' Preparation Theorem}, but
\textit{only with power series expansions}, and
then to apply it to prove the Approximation Problem, Cousin I/II Problems,
 $\delbar$-equation (for functions), 
holomorphic extensions (interpolations), and
 Levi's (Hartogs' Inverse)  Problem for unramified Riemann domains
 (over $\C^n$)\footnote{In the present, ``over $\C^n$'' will be
abbreviated, unless necessary.}
 (see Theorem \ref{thmchi} and \S\ref{levip});
 they constitute the main basic part of the theory
of several complex variables. 
The proofs are even simpler than those in the standard references
(cf., e.g., Gunning--Rossi \cite{guro}, Grauert--Remmert \cite{grr},
H\"ormander \cite{hor}, Noguchi \cite{nog16}).

Note that the present approach enables us to complete
 the proofs of those problems
in quite an elementary way
 without Weierstrass' Preparation Theorem or
 the cohomology theory\footnote{\, We use only the
$1$-cocyle class space
$H^1(\ast, \O)$ as a complex vector space.
} of Cartan--Serre, nor $L^2$-$\delbar$ method 
of H\"ormander.
The present paper came out from the study of the degree structure
of a generator system of a coherent analytic sheaf by \cite{nog15}.

Let $\Omega$ denote a domain of $\C^n$ with the structure sheaf
 $\O=\O_\Omega$.
For a holomorphic function $f \in \O(\Omega)$ in $\Omega$ we write
$\underline{f} \in \Gamma(\Omega, \O)$ for the induced
 sheaf-section of $\O$ and $\underline{f}_z$ for the germ of $f$
 at $z \in \Omega$.
Let $\sF$ be an analytic sheaf on $\Omega$ (i.e., a sheaf of $\O$-modules
over $\Omega$), and let
$\xi_j \in \Gamma(\Omega, \sF)$, $1 \leq j \leq q$, be finitely many
sections on $\Omega$. Then
the relation sheaf $\sR(\xi_1, \ldots, \xi_q)$ of
$\{\xi_j\}_{j=1}^q$ is a subsheaf of $\O^q$ consisting
of those germ-vectors
 $(\underline{f_1}_z, \ldots, \underline{f_q}_z)\in \O_z^q$
such that
\begin{equation}\label{rel}
 \underline{f_1}_z \xi_1(z)+ \cdots + \underline{f_q}_z \xi_q(z)=0,
\quad z \in \Omega.
\end{equation}

Now we formulate:
\begin{thm}[Weak Coherence]\label{weakc}
Let $S \subset \Omega$ be a complex submanifold.\footnote{\,A
 complex submanifold is not necessarily connected in this paper.}%
\begin{enumerate}
\setlength{\itemsep}{-3pt}%
\item
The non-singular geometric ideal sheaf $\sI\ang{S}$ is locally finite.
\item
Let $\left\{\underline{\sigma_j} \in \Gamma(\Omega, \sI \ang{S}
) : 1 \leq j \leq N\right\}$ be a finite generator system of
$\sI\ang{S}$ on $\Omega$ with $\sigma_j \in \O(\Omega)$:
 i.e.,
\[
 \sI\ang{S}  = \sum_{j=1}^N \O \cdot
     \underline{\sigma_j} .
\]
Then, the relation sheaf
$\sR(\underline{\sigma_1}, \ldots, \underline{\sigma_N})$
is locally finite.
\end{enumerate}
\end{thm}

We give a proof of this theorem in \S2. In \S3 we will apply it
to prove Oka's J\^oku-Ik\^o\footnote{\label{jkik}\,This
 is a method or a principle
of K. Oka all through his series of papers \cite{ok1}---\cite{ok9}
such that to solve a problem on a difficult domain one
embeds the domain into a higher dimensional polydisk, 
extends the problem on the polydisk, and then solves it by making use
of the simple shape of the polydisk (cf.\ \cite{nog16}).} (see
Lemma~\ref{ojki}),
 and then we will give a unified proof
for Cousin I/II Problems, and $\delbar$-equation for functions
in \S4 (Theorem \ref{thmchi}) by combining
 the Weak Coherence Theorem \ref{weakc}
 with a method of {\em cuboid induction on dimension}; these
yield $H^1(\Omega, \O)=0$ 
 for a holomorphically convex domain $\Omega$ (Lemma~\ref{h1}),
 which suffices to derive Oka's Heftungslemma or Grauert's
finiteness theorem for $\O$
on a strongly pseudoconvex domain (Theorem \ref{gr}).
In \S\ref{levip} we finally give  the solution of Levi's (Hartogs' Inverse) Problem
for unramified Riemann domains.

\section{Proof of Theorem \ref{weakc}}
(i) We take an arbitrary point $a \in \Omega$.

Case of $a \not\in S$: Since $S$ is closed, 
there is a neighborhood $U \subset \Omega$ of $a$ with
$U \cap S=\emptyset$. Then,
\[
 \sI\ang{S}_x=\O_x=1 \cdot \O_x, \quad \fa x \in U,
\]
and therefore, $\{1\}$ is a finite generator system of
$ \sI\ang{S}_x$ on $U$.

Case of $a \in S$: There is a holomorphic local coordinate neighborhood
$U$ of $a$ with $z=(z_1, \ldots, z_n)$ such that
\begin{align}\label{localcord}
a &= (0, \ldots, 0)\in U=\PD(0; (r_j)), \\ \notag
 S \cap U &=\{z=(z_j) \in U: z_1=\cdots =z_q=0\} \quad (1 \leq \ex q \leq n),
\end{align}
where $\PD(0; (r_j))=\{(z_j)\in \C^n: |z_j|<r_j, 1 \leq j \leq n\}$
is a polydisk with center at $0$.
Let $\underline{f}_b \in \sI\ang{S}_b$
($b \in U \cap S$) be any element.
With the coordinate system $(z_j)$ we write
 $b=(b_j)=(0, \ldots, 0, b_{q+1}, \ldots, b_n)$.
The function $f$ is represented by a unique
power series expansion, $f(z)=\sum_{\nu \in \Z_+^n} c_\nu(z-b)^\nu$,
which decomposes to
\begin{align*}
f(z) &= \sum_{\nu=(\nu_1, \nu') \in \Z_+^n, \nu_1>0} c_\nu(z-b)^\nu+
\sum_{\nu=(\nu_1, \nu') \in \Z_+^n, \nu_1=0} c_\nu(z-b)^\nu\\
&= \left(\sum_{\nu=(\nu_1,\nu') \in \Z_+^n, \nu_1>0} c_\nu z_1^{\nu_1-1}
(z'-b')^{\nu'}\right) z_1+
\sum_{\nu' \in \Z_+^{n-1}} c_{0\nu'} (z'-b')^{\nu'}.
\end{align*}
Here we put $\nu'=(\nu_2, \ldots, \nu_n), z'=(z_2, \ldots, z_n)$, and
$b'=(b_2, \ldots, b_n)$.
Setting
\begin{align*}
h_1(z_1,z') &=\left(\sum_{\nu=(\nu_1,\nu') \in \Z_+^n, \nu_1>0}
 c_\nu z_1^{\nu_1-1}(z'-b')^{\nu'}\right), \\
g_1(z')&=
\sum_{\nu' \in \Z_+^{n-1}} c_{0\nu'} (z'-b')^{\nu'},
\end{align*}
we have
\begin{equation}\label{div1}
 f(z_1, z')={h_1}(z_1, z') \cdot z_1+{g_1(z')}.
\end{equation}
For $g_1(z')$ we apply a similar decomposition with respect to
variable $z_2$, so that
\begin{equation*}
{g_1}(z')={h_2} \cdot z_2+ {g_2}(z''), \quad z''=(z_3, \ldots, z_n).
\end{equation*}
Repeating this process, we get
\[
f(z)= \sum_{j=1}^q h_j(z)\cdot z_j +g_q(z_{q+1}, \ldots, z_n).
\]
If $z_1=\cdots=z_q=0$, then $f(z)=0$, and so
$g_q(z_{q+1}, \ldots, z_n) =0$. Therefore,
\[
f(z)= \sum_{j=1}^q h_j(z)\cdot z_j.
\]
Thus, 
\begin{equation}\label{qgen}
\sI\ang{S}|_U= \sum_{j=1}^q \O_U \cdot \underline{z_j}.
\end{equation}

\medskip
(ii) We begin with the following lemma:
\begin{lem}\label{1-rel}
With the natural complex coordinate system
$z=(z_1, \ldots, z_n) \in \C^n$ we consider a relation sheaf
$\sR_p$  $(1 \leq p \leq n)$ defined by
\begin{equation}\label{s-rel}
\underline{f_1}_z \underline{z_1}_z +\cdots+
\underline{f_p}_z \underline{z_p}_z =0, \quad 
\underline{f_j}_z \in \O_{z}.
\end{equation}
Then $\sR_p$ is finitely generated on $\C^n$ by
\begin{equation}\label{trivialsol}
{T_{ij}}=
(0, \ldots,0, \overset{i\hbox{\scriptsize{\rm -th}}}{-\underline{z_j}}, 0,
\ldots, 0, \overset{j\mbox{\scriptsize{\rm -th}}}{\underline{z_i}},0,
\ldots, 0), \quad
1 \leq i<j \leq p.
\end{equation}
\end{lem}

We call ${T_{ij}}$ ($1 \leq i<j \leq p$) of \eqref{trivialsol} the
\textit{trivial solutions} of \eqref{s-rel} or of $\sR_p$.
In the case of $p=1$, we set the trivial solution to be $0$
as a convention.

\medskip
{\it Proof of Lemma \ref{1-rel}}\,: We use induction on $p \geq 1$.
The case of $p=1$ is clear.

Assuming that the case of $p-1 ~ (p \geq 2)$ holds,
we consider the case of $p$.  Set
$$
\Sig=\{(z_1, \ldots, z_n): z_1= \cdots =z_p=0\},
$$
and let $a \in \C^n$ be an arbitrary point.
If $a=(a_j) \not\in \Sig$, there is an $a_j \not=0$ ($1 \leq j \leq p$),
to say, $a_1 \not=0$. In a neighborhood $V$ of $a$,
 $z_1\not=0$. Then, \eqref{s-rel} is solvable with respect to
 $\underline{f_1}_z$:
\[
 \underline{f_1}_{z}= - \underline{f_2}_{z} \cdot
\frac{\underline{z_2}_{z}}{\underline{z_1}_z}- \cdots -
\underline{f_p}_{z} \cdot \frac{\underline{z_p}_{z}}{\underline{z_1}_z},
\quad \fa \underline{f_j}_z \in \O_{z}~ (2 \leq j \leq p),
 ~  z \in V.
\]
It follows that with $z \in V$,
\begin{align}\label{trivial1}
\left(\underline{f_j}_z\right)&=\left(- \sum_{j=2}^p \underline{f_j}_z \cdot
\frac{\underline{z_j}_{z}}{\underline{z_1}_z}, \underline{f_2}_z,
\ldots , \underline{f_p}_z\right)\\ \notag
&=  \sum_{j=2}^p  \frac{\underline{f_j}_z}{\underline{z_1}_z} \cdot
\left( - {\underline{z_j}_z}, 0, \ldots, 0 ,
 \overset{j\mbox{\scriptsize{\rm -th}}}{\underline{z_1}_z}, 0,
\ldots, 0\right)\\ \notag
&=  \sum_{j=2}^p - \frac{\underline{f_j}_z}{\underline{z_1}_z}
\cdot {T_{1j}}(z) \in \sum_{j=2}^p \O_{z} \cdot
 {T_{1j}}(z) .
\end{align}
Therefore, $\sR_p$ is generated by the trivial solutions
$\{{T_{1j}}\}_{2 \leq j \leq  p}$ on $V$.

If $a \in \Sig$, we decompose an element 
$\left( \underline{f_j}_a \right) \in \sR_{p\,a}$ in a polydisk neighborhood
$U$ of $a$ as in \eqref{div1}:
\[
  f_j(z_1, z')={h_j}(z_1,z')  z_1+{g_j(z')},
\quad z'=(z_2, \ldots, z_n), ~ 1 \leq j \leq p.
\]
For $z \in U$ one gets
\begin{align}\label{T1}
\left(\underline{f_j}_z\right) - \sum_{j=2}^p \underline{h_j}_z\,
T_{1j}(z) &=
\left(\underline{g_1}_z+ \sum_{j=1}^p \underline{h_j}_z\,
 \underline{z_j}_z,~ 
\underline{g_2}_z, \ldots, \underline{g_p}_z  \right) \\ \notag
&=\left(\underline{\tilde{g}_1}_z,
\underline{g_2}_z, \ldots, \underline{g_p}_z  \right).
\end{align}
Here, $\underline{\tilde{g}_1}_z=
\underline{g_1}_z+ \sum_{j=1}^p \underline{h_j}_z \underline{z_j}_z$.
Since $\left(\underline{\tilde{g}_1}_z,
\underline{g_2}_z, \ldots, \underline{g_p}_z  \right) \in {\sR_p}_z$,
\[
\underline{\tilde{g}_1}_z \,  \underline{z_1}_z+
\underline{g_2}_z \, \underline{z_2}_z+
 \cdots + \underline{g_p}_z   \underline{z_p}_z=0 .
\]
The second term and so forth of the right-hand side of the equation
above do not contain variable $z_1$, and so $\underline{\tilde{g}_1}_z=0$
is deduced. Thus,
\[
\underline{g_2}_z \,  \underline{z_2}_z+
 \cdots + \underline{g_p}_z \, \underline{z_p}_z=0.
\]
This is the case of $p-1$ after changing the indices of variables.
Therefore, the induction hypothesis implies that
$\left(0, \underline{g_2}_z, \ldots, \underline{g_p}_z\right)$
is represented as a linear sum of
$T_{ij}(z),~ 2 \leq i <j  \leq p$, with coefficients
in $\O_z$.
Combining this with \eqref{T1}, we see that
$\left(\underline{f_j}_z\right)$ is represented as a linear sum of
$T_{ij}(z),~ 1 \leq i<j  \leq p$, with coefficients
in $\O_z$.
\hfill$\triangle$

\medskip
\textit{Continued proof of} (ii):
Set
$\sR=\sR(\underline{\sigma_1}, \ldots, \underline{\sigma_N})$.
We consider the relation
\begin{equation}\label{sigma-rel}
\underline{f_1}_z\, \underline{\sigma_1}_z +\cdots+
\underline{f_N}_z \, \underline{\sigma_N}_z =0, \quad 
\underline{f_j}_z \in \O_z.
\end{equation}
We set the trivial solutions of this equation as follows:
\[
\tau_{ij}=
(~\ldots, \, \overset{i\mbox{\scriptsize{-th}}}{- \underline{\sigma_j}},
\ldots, \overset{j\mbox{\scriptsize{-th}}}{\underline{\sigma_i}},
\ldots ~), \quad
1 \leq i<j \leq N.
\]

We take an arbitrary point $a \in \Omega$.
If $a \not\in S$, then some $\sigma_j(a)\not=0$, to say,
$\sigma_1(a) \not=0$.  As in \eqref{trivial1}, one sees that
$\sR$ is generated by $\{\tau_{1j} \}_{j=2}^N$ on
a neighborhood of $a$.

If $a \in S$, we take a holomorphic local coordinate
system $z=(z_1, \ldots, z_n)$ in a polydisk neighborhood
 $\PD$ as in \eqref{localcord}:
\begin{align*}
a &= (0, \ldots, 0),\\
 S\cap \PD &= \{(z_1, \ldots, z_n) \in \PD: z_1= \cdots=z_q=0\}\quad
(1 \leq  \ex q \leq n).
\end{align*}
It follows from \eqref{qgen} and the assumption that
\[
\sI\ang{S}|_\PD= 
\sum_{j=1}^q \O_\PD \cdot \underline{z_j}=
\sum_{j=1}^N \O_\PD \cdot \underline{\sigma_j}|_\PD .
\]
Thus, we may assume without loss of generality that
\begin{align*}
\sigma_j &= z_j, \quad 1 \leq j \leq q~ (\hbox{on } \PD),\\
\sigma_i &=\sum_{j=1}^q a_{ij} \, z_j, \quad a_{ij}\in
 \O(\PD),~ q+1 \leq i \leq N ~ (\hbox{on } \PD).
\end{align*}
Set
\begin{equation}\label{non-tri}
\phi_i=(-\underline{a_{i1}}, \ldots, -\underline{a_{iq}}, 0,
 \ldots, 0,
 \overset{i\mbox{\scriptsize{-th}}}{{1}}, 0,  \ldots , 0 )
\in \Gamma(\PD, \sR), \quad q+1 \leq i  \leq N.
\end{equation}

We deduce from \eqref{sigma-rel} with $z \in \PD$ that
\begin{equation}\label{q-rel}
 \left(\underline{f_1}_z+\sum_{i=q+1}^N \underline{f_i}_z\,
\underline{a_{i1}}_z\right)\underline{z_1}_z+ \cdots +
 \left(\underline{f_q}_z+\sum_{i=q+1}^N \underline{f_i}_z\,
\underline{a_{iq}}_z\right)\underline{z_q}_z=0.
\end{equation}
By Lemma \ref{1-rel},
\[
 \left(\underline{f_1}_z+\sum_{i=q+1}^N \underline{f_i}_z\,
\underline{a_{i1}}_z , \ldots,
\underline{f_q}_z+\sum_{i=q+1}^N 
 \underline{f_i}_z \, \underline{a_{iq}}_z, 0, \ldots, 0\right)
\]
is a linear sum of $\tau_{jk}(z)$, $1 \leq j < k \leq q$,
with coefficients in $\O_z$.
Therefore there are
$\underline{b_{jk}}_z \in \O_z$, $1 \leq j < k \leq q$, such that
\begin{align}\label{q-sol}
\sum_{1 \leq j<k\leq q} \underline{b_{jk}}_z \tau_{jk}(z)
&=
\left(\underline{f_1}_z+\sum_{i=q+1}^N \underline{f_i}_z\,
\underline{a_{i1}}_z , \ldots,
\underline{f_q}_z+\sum_{i=q+1}^N 
 \underline{f_i}_z \, \underline{a_{iq}}_z, 0, \ldots, 0\right) \\
 \notag
&=\left(\underline{f_1}_z , \ldots, \underline{f_q}_z, 0, \ldots, 0 \right)
+ \sum_{i=q+1}^N \underline{f_i}_z \, \left(
\underline{a_{i1}}_z, \ldots,\underline{a_{iq}}_z, 0, \ldots, 0\right) .
\end{align}
By making use of \eqref{non-tri} we get
\begin{align}\label{q-sol2}
\left(\underline{f_1}_z, \ldots, 
\underline{f_q}_z, \ldots, \underline{f_N}_z\right)
= \sum_{1 \leq j<k\leq q} \underline{b_{jk}}_z \,
 \tau_{jk}(z)
+ \sum_{i=q+1}^N \underline{f_i}_z\, \phi_i(z) .
\end{align}
Thus, $\sR$ is generated on $\PD$ by
\begin{equation}\label{1syzy}
 \tau_{jk}, ~  \phi_i, \quad
 1 \leq j < k \leq q, ~ q+1 \leq i \leq N.
\end{equation}
This finishes the proof. \qed

\begin{rmk}\rm\begin{enumerate}
\setlength{\itemsep}{-3pt}
\item
In the Weak Coherence Theorem~\ref{weakc} it is a point
to assume that $\{\underline{\sigma_j}\}_{j=1}^N$ is a generator
system of $\sI\ang{S}$; otherwise, the proof above does not
     work even if $S$ is non-singular.
\item
It is an advantage of the above method to the general First Coherence
Theorem of Oka that we have an explicit system of generators \eqref{1syzy}.
\end{enumerate}
\end{rmk}

\section{Oka's J\^oku-Ik\^o}
The term ``J\^oku-Ik\^o'' was used by K. Oka since
he wrote the first paper of the series in 1936
and retained this principle all through
his works (\cite{ok1}---\cite{ok9})
(see footnote \ref{jkik});
The aim of the present section is to prove
Oka's J\^oku-Ik\^o, Lemma~\ref{ojki} below only by making use of
Theorem~\ref{weakc} combined with Cousin's integral \eqref{cous}.
The technics may be essentially similar to those
in some references, e.g., Nishino \cite{nis} and Noguchi \cite{nog16},
but they are not in a suitable form for our purpose.

\subsection{Syzygy for non-singular geometric ideal sheaves}
We begin with:
\begin{defn}
\label{closedcuboid}\rm
A cuboid $E$ is a bounded open or closed subset of $\C^n$
with the boundary
parallel to the real and imaginary axes of $z=(z_1, \ldots, z_n) \in
 \C^n$.
In the case of $n=1$, $E$ is called a rectangle.
When $E$ is a closed cuboid,
 we allow the widths of some edges to degenerate to $0$,
and call the number of edges of $E$ of positive widths the dimension of
$E$, denoted by $\dim E$.
\end{defn}

Let $\Omega \subset \C^n=\C^{n-1} \times \C$ be a domain and
let $E', E'' \Subset \Omega$\,\footnote{ The symbol ``$\Subset$''
stands for that the inclusion is relatively compact.}
 be two closed cuboids as follows:
There are a closed cuboid
$F \Subset \C^{n-1}$ and two adjacent closed rectangles
$E_n', E_n'' \Subset \C$ sharing a side $\ell$, and
\begin{align}\label{adjc}
E' =F \times E_n',\quad  E''=F \times E_n'' , \quad
\ell =E_n' \cap E_n'' .
\end{align}
\begin{figure}[h] 
\begin{center}
\includegraphics{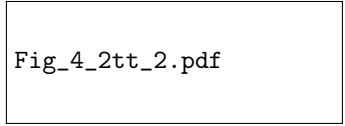}
\end{center}
\caption{Adjacent closed cuboids} 
\label{cuboid1}
\end{figure}

We now recall:
\begin{lem}[Cartan's Merging Lemma]
\label{cartan-merge}
Let $E', E''\Subset \Omega$ be adjacent
closed cuboids as in \eqref{adjc}, and let
$\sF$ be an analytic sheaf on $\Omega$.
Let $\{\sigma'_j \in \Gamma(U', \sF): 1 \leq j \leq p'\}$
(resp.\ $\{\sigma''_k \in \Gamma(U'', \sF),
 1 \leq k \leq p''\}$)
be a finite generator system of $\sF$ on $E'$ (resp.\ $E''$).\footnote{
This means that they are defined so in some neighborhoods of $E'$ and $E''$,
 respectively; this expression is the same through the paper.}

Moreover, assume that
there are holomorphic functions
$a_{jk}, b_{kj} \in \O(E' \cap E'')$, $1 \leq j \leq p'$,
$1 \leq k \leq p''$, such that
$$
\sigma'_j=\sum_{k=1}^{p''} \underline{a_{jk}} \cdot
 \sigma''_k, \quad
\sigma''_k=\sum_{j=1}^{p'} \underline{b_{kj}} \cdot
 \sigma'_j
\quad (\hbox{on } E' \cap E'').
$$

Then, there exists a merged finite generator system
$\{\sigma_l \in \Gamma(E'\cup E'', \sF):
 1 \leq l \leq p'+p''\}$ on $E' \cup E''$.
\end{lem}

This is due to H. Cartan 1940;
a rather simplified proof of it
can be found in \cite{nog16}, ``Added at galley-proof''.

\begin{lem}[Oka's Syzygy]
\label{oka-syzygy}
Let $E \Subset \C^n$ be a closed cuboid.
\begin{enumerate}\setlength{\itemsep}{-3pt}
\item
Every locally finite analytic sheaf $\sF$ defined on $E$
(i.e., in a neighborhood of $E$)
 has a finite generator system on $E$.
\item
Let $\sF$ be an analytic sheaf on $E$ with a finite generator system
$\{\sigma_j\}_{1 \leq j \leq N}$ on $E$ such that the relation sheaf
$\sR(\sigma_1, \ldots, \sigma_N)$ is locally finite.

Then  for every section $\sigma \in \Gamma(E, \sF)$
there are holomorphic functions $a_j \in \O(E)$, $1 \leq j \leq N$,
such that
\end{enumerate}\vspace{-10pt}
\begin{equation}\label{repres}
\sigma= \sum_{j=1}^N ~\underline{a_j}
 \cdot \sigma_j \quad (\hbox{on } E).
\end{equation}
\end{lem}
\begin{pf}
The proof is carried out in the same way as in 
\cite{nis}, or \cite{nog16} Lemma~4.3.7
except for the use of the vanishing $H^1(U, \O)=0$  for an affine convex cylinder
domain $U \subset \C^n$, which we replace by Cousin's integral
\eqref{cous} as follows.
Suppose that $E$ is a closed cuboid such that
\begin{align}\label{cuboid}
  E = F \times \{z_n: | \re z_n| \leq T,~ |\im z_n|\leq \theta\},\qquad
T>0, ~  \theta \geq 0.
\end{align}
Set $E_0= F \times \{z_n: \re z_n=0,\, |\im z_n|\leq \theta\}$,
and let $\varphi(z', z_n) \in \O(E_0)$. Then there is a small
$\delta>0$ such that $\varphi(z', z_n)$ is defined on
\[
 F \times \{z_n: |\re z_n| \leq \delta,~ |\im z_n|\leq \theta+\delta\}.
\]
Set
\begin{align*}
\ell &= \{z_n: \re z_n=0, ~  -\theta - \delta \leq \im z_n\leq
 \theta+\delta\}, \\
E_1 &= F \times \{z_n: -T \leq \re z_n \leq \delta,~ |\im z_n|\leq \theta\},\\
E_2 &= F \times \{z_n: -\delta \leq \re z_n \leq T,~ |\im z_n|\leq
 \theta\},
\end{align*}
where $\ell$ is positively oriented as $\im z_n$ increases.
We define Cousin's integral of $\varphi(z',z_n)$ along $\ell$ by
\begin{equation}\label{cous}
 \Phi(z', z_n)= \frac{1}{2\pi i}
\int_\ell ~ \frac{\varphi(z', \zeta_n)}{\zeta_n -z_n}~
d\zeta_n.
\end{equation}
Then $\Phi(z', z_n)$ is holomorphic on
 $(E_1 \cup E_2) \setminus (F \times \ell)$.
After analytic continuations we obtain $\Phi_j(z',z_n) \in \O(E_j)$
($j=1,2$) satisfying
\begin{equation}\label{cousd}
\Phi_1(z', z_n) - \Phi_2(z', z_n) =\varphi(z', z_n), \quad
 (z',z_n) \in E_1 \cap E_2.
\end{equation}
We call this the {\it Cousin decomposition} of $\varphi(z',z_n)$.

The rest is the same as in the proof of \cite{nog16} Lemma~4.3.7.
\end{pf}

By the Weak Coherence Theorem~\ref{weakc} and Lemma~\ref{oka-syzygy}
we have:
\begin{thm}[Syzygy for $\sI\ang{S}$]\label{gsyzygy}
Let $S$ be a complex submanifold of a neighborhood of
a closed cuboid $E~(\Subset \C^n)$.
\begin{enumerate}\setlength{\itemsep}{-3pt}
\item
$\sI\ang{S}$ has a finite generator system on $E$.
\item
Let $\{\underline{\sigma_j}\}_{1 \leq j \leq N}$ be a finite
generator system of $\sI\ang{S}$ on $E$
with $ \sigma_j \in \O(E)$.
Then for every
$\underline{\sigma} \in \Gamma(E, \sI\ang{S})$
($\sigma \in \O(E)$) there are holomorphic functions
$a_j \in \O(E)$, $1 \leq j \leq N$, such that
\end{enumerate}\vspace{-10pt}
\begin{equation}\label{geom-repres}
{\sigma}= \sum_{j=1}^N ~{a_j}
 \cdot {\sigma_j} \quad (on ~E).
\end{equation}
\end{thm}

\subsection{Oka's J\^oku-Ik\^o}

Let $\rP$ be an open cuboid in $\C^n$, and let $S \subset \rP$
be a complex submanifold. The following is fundamental in the Oka theory.
\begin{lem}[Oka's J\^oku-Ik\^o]\label{ojki}
Let $E \Subset \rP$ be a closed cuboid.
Then for every holomorphic function
$g$ on $E \cap S ~(\Subset S)$\footnote{\label{nbd}%
With this writing we mean that
$g$ is a holomorphic function in a neighborhood $V$ of
 $E \cap S$ in $S$. The notation will be used in sequel.}
 there exists an element $G \in \O(E)$  satisfying
$$
G|_{E \cap S}=g|_{E \cap S}.
$$
Here, the equality holds in a neighborhood of $E \cap S$
in $S$.\footnote{
The formulation of this lemma and the proof below should be new.
}
\end{lem}

{\em Remark.} We call $G$ above a {\em solution} on $E$.
\begin{pf}  Notice that in the case of $E \cap S=\emptyset$, $G$ can be any
 holomorphic function on $E$, and the statement is true.
We use induction on $\dim E$.

(a) Case of $\dim E=0$: Since $E$ consists of one point,
the assertion is clear.

(b) Case of $\dim E=\nu$ ($\nu \geq 1$) with the induction hypothesis
that the case of $\dim E=\nu -1$ is true: 
 By Theorem \ref{gsyzygy}~(i) there is
a finite generator system 
$\{\underline{\sigma_j}\}_{j=1}^N$ of $\sI\ang{S}$
on a neighborhood $W (\subset \rP)$ of $E$ with $\sigma_j \in \O(W)$.

We may assume that $E$ is taken as in \eqref{cuboid}.
We set
\begin{equation}\label{cuboidf}
 E_t=\{z=(z', z_n) \in E: \re z_n=t\}, \quad -T \leq t \leq T.
\end{equation}
Since $\dim E_t=\nu-1$, the induction hypothesis implies that
there is a solution $G_t \in \O(E_t)$ satisfying
$G_t|_{S \cap E_t}=g|_{S \cap E_t}$.
By the Heine--Borel Theorem there is a finite partition
\begin{align}\label{part0}
-T &=t_0 < t_1 < \cdots < t_L=T, \\ \notag
E_\alpha &:= \{z=(z', z_n) \in E: 
 t_{\alpha-1} \leq  \re z_n \leq t_\alpha \}, \quad
1 \leq \alpha \leq L,
\end{align}
such that there are solutions $G_\alpha\in \O(E_\alpha)$ satisfying
$$
 G_\alpha|_{S \cap E_\alpha}=g|_{S \cap E_\alpha}.
$$
Therefore, $\underline{G_{\alpha+1} - G_\alpha} \in
\Gamma(E_\alpha \cap E_{\alpha+1}, \sI\ang{S})$.
It follows from Theorem \ref{gsyzygy}~(ii) that there are
$a_{\alpha j} \in \O(E_\alpha \cap E_{\alpha+1})$ ($1 \leq j \leq N$)
satisfying
\begin{equation}\label{1cyc}
G_{\alpha+1}-G_{\alpha}= \sum_{j=1}^N a_{\alpha j} \sigma_j \quad
(\hbox{on }  E_\alpha \cap E_{\alpha+1}).
\end{equation}
By the Cousin decomposition \eqref{cous} of $a_{\alpha j}$
we write
\begin{equation}\label{c-decom}
a_{\alpha j}= b_{\alpha j} - b_{\alpha+1 j} \quad
 (\hbox{on } E_\alpha \cap E_{\alpha+1}), \quad
b_{\alpha j} \in \O(E_\alpha),~b_{\alpha+1 j} \in \O(E_{\alpha+1}).
\end{equation}
Then,
\begin{equation}\label{conec}
G_{\alpha}+ \sum_{j=1}^N b_{\alpha j}  \sigma_j=
G_{\alpha+1} + \sum_{j=1}^N b_{\alpha+1 j} \sigma_j \quad
(\hbox{on } E_\alpha \cap E_{\alpha+1}).
\end{equation}
Thus this yields a solution $H_{\alpha+1}$ on
$E_{\alpha} \cup E_{\alpha+1}$;  for this procedure we say
that we merge the solutions $G_\alpha$ and $G_{\alpha+1}$
to obtain a solution $H_{\alpha+1}$ on $E_\alpha \cup E_{\alpha+1}$.

Starting from $\alpha=1$, we merge $G_1$ and $G_2$ to obtain
a solution $H_2$ on $E_1 \cup E_2$. We then merge
$H_2$ and $G_3$ to obtain a solution $H_3$ on $E_1\cup E_2 \cup E_3$.
Repeating this procedure up to $\alpha=L-1$, we obtain a solution
$H_L$ on  $E=\bigcup_{\alpha=1}^L E_\alpha$, and set $G=H_L$:
This finishes the proof of Lemma~\ref{ojki}.
\end{pf}

\begin{rmk}\label{cdi}\rm
We call the above induction argument  {\it cuboid induction on dimension},
which will be used furthermore in the sequel.
\end{rmk}

It is well-known that Oka's J\^oku-Ik\^o Lemma~\ref{ojki}
immediately implies (cf., e.g., \cite{nog16} Lemma~4.4.17):

\begin{thm}[Runge--Weil--Oka Approximation]\label{rwo}
Let $\Delta \Subset \Omega$ be an analytic polyhedron
of a domain $\Omega~(\subset \C^n)$.
Then every holomorphic function on the closure $\bar\Delta$ is uniformly
 approximated on $\bar\Delta$ by elements of $\O(\Omega)$.
\end{thm}

\section{Cousin I/II, $\delbar$, Extension and Levi's (Hartogs' Inverse) Problems}
The aim of this section is to show how the result obtained
in the previous section is applied to solve the titled problems.

\subsection{Cousin I/II, and $\delbar$-Equation
}\label{uni}
We will give one unified proof to all of the three problems.
We recall them: Let $\Omega \subset \C^n$ be a domain, let
 $\Omega=\bigcup_{\alpha \in \Lambda} U_\alpha$ be an open covering, and let $\sM(U_\alpha)$
denote the set of all meromorphic functions in $U_\alpha$.
\begin{enumerate}\setlength{\itemsep}{-3pt}
\item[I]
({\bf Cousin I}) For given $f_\alpha \in \sM(U_\alpha)$ ($\alpha \in \Lambda$)
satisfying $f_\alpha - f_\beta \in \O(U_\alpha \cap U_\beta)$ (Cousin-I data),
find $F \in \sM(\Omega)$ (called a {\it solution}) with
 $F|_{U_\alpha}- f_\alpha \in \O(U_\alpha)$
for all $\alpha \in \Lambda$.
\item[II]
({\bf Cousin II, Oka Principle}) Here we assume that $U_\alpha$ are simply-connected.
Let  $f_\alpha \in \sM^*(U_\alpha)$ ($\alpha \in \Lambda$) be
 locally non-zero meromorphic functions
satisfying
\begin{enumerate}\setlength{\itemsep}{-3pt}
\item
$f_\alpha / f_\beta \in \O^*(U_\alpha \cap U_\beta)$ (nowhere
	 vanishing holomorphic functions) (Cousin-II data),
\item
(Topological condition) there are nowhere vanishing continuous functions
$\psi_\alpha \in \sC^*(U_\alpha)$ with
$\psi_\alpha/\psi_\beta=f_\beta/f_\alpha$ on $U_\alpha \cap U_\beta$.
\end{enumerate}
Find $F \in \sM^*(\Omega)$ with $F|_{U_\alpha}/ f_\alpha \in \O^*(U_\alpha)$
for all $\alpha \in \Lambda$.

We may take a continuous branch $h_\alpha:=\log \psi_\alpha$
in each $U_\alpha$. It follows that
\end{enumerate}
\begin{equation}\label{cocy}
{h_\alpha}|_{U_\alpha \cap U_\beta} -
{h_\beta}|_{U_\alpha \cap U_\beta} \in \O(U_\alpha \cap U_\beta).
\end{equation}
\begin{enumerate}\setlength{\itemsep}{-3pt}
\item[{}]
The problem to find the solution $F$ above is reduced to
the following problem:

$(\star)$  For a given family $\{h_\alpha\}$ of continuous functions
satisfying \eqref{cocy}, find a continuous function
$H$ (called a solution) in $\Omega$ such that for every $U_\alpha$
\[
 H|_{U_\alpha} - h_\alpha \in \O(U_\alpha).
\]
\item[III]
({\bf $\delbar$-Equation}) For a given $C^\infty$-$(0,1)$-form $u$ on $\Omega$
 with $\delbar u=0$, find a $C^\infty$-function $g$
(called a {\em solution}) on $\Omega$ with
$\delbar g=u$.

 Locally, by Dolbeault's lemma, there is a solution $f$
of this problem in a neighborhood of a point of $\Omega$.
Thus, there are an open covering $\{U_\alpha\}_{\alpha\in\Lambda}$ of
$\Omega$ and $C^\infty$-functions $g_\alpha$ on $U_\alpha$ such that
$\delbar g_\alpha=u|_{U_\alpha}$. Then, the present problem is
equivalent to find a $C^\infty$-function $G$ ({\it solution})
on $\Omega$ with $G|_{U_\alpha} - g_\alpha \in \O(U_\alpha)$ for all
$\alpha \in \Lambda$.
\end{enumerate}

{\bf Convention.}  For a unified treatment for the above problems,
 we introduce
an ``{\it argument} $\chi$'' representing one of I---III above,:
Problem-$\chi$ means one of Problems I---III above, 
where Problem-II means Problem ($\star$),
and a $\chi$-solution
means a {\it solution} of the corresponding Problem-$\chi$.

\begin{rmk}\rm
If $\Psi$ is so obtained in Cousin-II Problem above, then
$F_1=f_\alpha e^{\log \psi_\alpha -\Psi} \in \sM^*(\Omega)$ satisfies
the required property for $F$. Then we have a homotopy,
\[
 F_t=f_\alpha e^{\log \psi_\alpha  -t\Psi}, \quad 0 \leq t \leq 1,
\]
from the topologically assumed function
 $F_0(=f_\alpha \psi_\alpha)$
to an aimed analytic (meromorphic) function $F_1$.
\end{rmk}

\begin{rmk}\label{comm}\rm
The common property of Problem-$\chi$ that we will use is the following:
If $f$ and $f'$ are  two solutions of Problem-$\chi$ on an open set $U$
in general, then $f-f' \in \O(U)$.
\end{rmk}

We begin with:
\begin{lem}\label{chilem}
Let $\rP$ be an open cuboid in $\C^n$ and let $S$ be a complex submanifold
of $\rP$. We consider Problem-$\chi$ defined on $S$. 
 Let $E \Subset \rP$ be
a closed cuboid. Then there is a $\chi$-solution on
 $E \cap S (\Subset S)$.\footnote{ Cf.\ footnote \ref{nbd} at
p.~\pageref{nbd}.}
\end{lem}
\begin{pf}
We use cuboid induction on dimension.

(a) Case of $\dim E=0$: It is clear by definition.

(b) Case of $\dim E=\nu (\nu \geq 1)$ with the induction hypothesis that the case
of $\dim E=\nu-1$ holds:
Without loss of generality we may assume that $E$ is given as in
\eqref{cuboid}, and let $E_t$ be as in \eqref{cuboidf}.
Since $\dim E_t=\nu-1$, the induction hypothesis implies the existence
of a $\chi$-solution $\Phi_t$ on $E_t \cap S~(\Subset S)$.
 Then, by the Heine-Borel
Theorem there are a partition of $[-T, T]$, $E_\alpha$ ($1 \leq \alpha \leq L$)
 as in \eqref{part0}, and $\chi$-solutions
 $\Phi_\alpha$ on $E_\alpha \cap S (\Subset S)$.

If $E_\alpha \cap E_{\alpha+1}\cap S \not=\emptyset$, we say that
$E_\alpha$ and $E_{\alpha+1}$ is pairwise connected on $S$.
It is sufficient to prove the existence of a $\chi$-solution for
each maximal sequence of $E_\alpha$ pairwise connected on $S$,
\begin{equation}\label{mcs}
 E_{\alpha_0} \cup E_{\alpha_0 +1} \cup \cdots \cup E_{\alpha_1}.
\end{equation}
For simplicity we suppose that $\alpha_0=1$.
It follows from  Remark~\ref{comm} that for $ 1 \leq \alpha \leq \alpha_1$
\begin{equation}\label{diff}
\Phi_{\alpha+1} - \Phi_\alpha \in
\Gamma(E_\alpha \cap E_{\alpha+1}\cap S, \O_S).
\end{equation}
By Oka's J\^oku-Ik\^o Lemma~\ref{ojki},
 there is a holomorphic function
 $H_\alpha \in \O(E_\alpha \cap E_{\alpha+1})$ such that
\begin{equation}\label{jki}
H_\alpha|_{E_\alpha \cap E_{\alpha+1} \cap S}=
\Phi_{\alpha+1} - \Phi_\alpha. 
\end{equation}
By the Cousin decomposition of $H_\alpha$ as in \eqref{cous} we have
$\tilde{H}_\alpha \in \O(E_\alpha)$ and
$\tilde{H}_{\alpha+1} \in \O(E_{\alpha+1})$ such that
\begin{equation}\label{dec}
 H_\alpha = \tilde{H}_\alpha - \tilde{H}_{\alpha+1} \quad
(\hbox{on } E_\alpha \cap E_{\alpha+1}).
\end{equation}

We infer from \eqref{dec} and \eqref{ext} that
\begin{equation}\label{merg}
\Phi_\alpha+\tilde{H}_\alpha|_{E_\alpha\cap S} =
\Phi_{\alpha+1}+\tilde{H}_{\alpha+1}|_{E_{\alpha+1}\cap S}
\quad \hbox{on } E_\alpha \cap E_{\alpha+1} \cap S \, (\Subset S).
\end{equation}
Note that $\Phi_\alpha+\tilde{H}_\alpha|_{E_\alpha\cap S}$ 
(resp.\ $\Phi_{\alpha+1}+\tilde{H}_{\alpha+1}|_{E_{\alpha+1}\cap S}$)
is a $\chi$-solution on $E_\alpha \cap S (\Subset S)$
(resp. $E_{\alpha+1} \cap S (\Subset S)$).
Thus, from \eqref{merg} we obtain a merged $\chi$-solution
 $\Psi_{\alpha+1}$ on
$(E_\alpha\cup E_{\alpha+1}) \cap S (\Subset S)$ from $\Phi_\alpha$ and
$\Phi_{\alpha+1}$.

Now, from $\Phi_1$ and $\Phi_2$ we obtain a merged $\chi$-solution
$\Psi_2$ on $(E_1 \cup E_2)\cap S (\Subset S)$. We then obtain a merged
$\chi$-solution $\Psi_3$ on $(E_1 \cup E_2 \cup E_3)\cap S (\Subset S)$
from $\Psi_2$ and $\Phi_3$, and so on; we obtain a
$\chi$-solution on $(\bigcup_{\alpha=1}^{\alpha_1}E_\alpha)\cap S (\Subset S)$.
\end{pf}

\begin{defn}\label{hconv}\rm
A complex manifold $M$ is said to be {\em holomorphically convex}
if for every compact subset $K \Subset M$ the holomorphically
convex hull of $K$ defined by
\[
 \hat{K}_M = \{a \in M: |f(a)| \leq \max\{|f(x)|: x \in K\}
\]
is again compact.
\end{defn}

\begin{thm}\label{thmchi}
Let $\Omega$ be a holomorphically convex domain.
Then Problem-$\chi$ on $\Omega$ has a $\chi$-solution on $\Omega$.
\end{thm}
\begin{pf}
We take an increasing sequence of analytic polyhedra of $\Omega$,
\begin{equation}\label{aph}
 \Delta_1 \Subset \Delta_2 \Subset \Delta_3 \Subset \cdots,
 \qquad \bigcup_{\nu=1}^\infty \Delta_\nu=\Omega.
\end{equation}
For each $\nu$ we let $\phi_\nu: \bar\Delta_\nu \to \overline{\PD}_\nu$ be
the Oka map (a holomorphic proper embedding)
of $\bar\Delta_\nu$ into a closed polydisk $\overline{\PD}_\nu$,
which extends from a neighborhood $U_\nu$ of $\bar\Delta_\nu$ into
a polydisk, biholomorphic to an open cuboid
 $\rP_\nu\,(\Supset \overline\PD_\nu)$.
Then, the image $\phi_\nu(U_\nu)$ is 
 a complex submanifold of $\rP_\nu$.
We identify $U_\nu$ with the image $\phi_\nu(U_\nu)$.

By Lemma \ref{chilem} there is a $\chi$-solution $G_\nu$
on every $\bar\Delta_\nu$. Put $F_1=G_1$ on $\bar\Delta_1$.
Suppose that $\chi$-solutions $F_\nu$ on $\bar\Delta_\nu$,
 $1 \leq \nu \leq \mu$, are determined so that
\begin{equation}\label{maj}
 \|F_{\nu+1} - F_\nu\|_{\bar\Delta_\nu} <  \frac{1}{2^\nu},\quad
1 \leq \nu \leq \mu.
\end{equation}
Let $G_{\mu+1}$ be a $\chi$-solution on $\bar\Delta_{\mu+1}$.
Since $G_{\mu+1}|_{\bar\Delta_\mu} -F_\mu \in \O(\bar\Delta_\mu)$,
by Theorem~\ref{rwo} there is an element
 $h_{\mu+1} \in \O(\bar\Delta_{\mu+1})$
with
\[
 \|G_{\mu+1}|_{\bar\Delta_\mu} - F_\mu -h_{\mu+1}\|_{\bar\Delta_\mu}
< \frac{1}{2^{\mu+1}}.
\]
Setting $F_{\mu+1}=G_{\mu+1} -h_{\mu+1}$, we see that \eqref{maj}
holds up to $\mu+1$.
Inductively, we have $\chi$-solutions $F_\nu$ on $\bar\Delta_\nu$
satisfying \eqref{maj}, and the series
\[
 F=F_\mu + \sum_{\nu=\mu}^\infty (F_{\nu+1}- F_\nu)
\]
converges  locally uniformly and the limit
gives rise to a $\chi$-solution on $\Omega$.
\end{pf}

\begin{rmk}\rm
 As easily seen, the above proof of Theorem~\ref{thmchi}
 works on Stein manifolds, which is defined by:
\begin{defn}\label{stein}\rm
 A complex manifold $M$ is said to be {\em Stein} if
the following conditions are satisfied.
\vspace{-3mm}
\begin{enumerate}\setlength{\itemsep}{-3pt}
\item
$M$ satisfies the second countability axiom:
\item
({\em Holomorphic separation})
For every two distinct two point $a, b \in M$ there is
a holomorphic function $f \in \O(M)$ on $M$ with
$f(a) \not= f(b)$:
\item
For every point $a \in M$ there are holomorphic functions
$f_j \in \O(M), 1 \leq j \leq n=\dim_a M$, such that
$df_1(a) \wedge \cdots \wedge df_n (a) \not=0$:
\item
$M$ is holomorphically convex (see Definition \ref{hconv}).
\end{enumerate}
\end{defn}
\end{rmk}

\subsection{Extension Problem}
By means of the Weak Coherence Theorem~\ref{weakc}
we consider the extension problem (interpolation problem)
from a complex submanifold in a holomorphically convex domain.
\begin{thm}\label{ext}
Let $\Omega \subset \C^n$ be a holomorphically convex domain and let
$S \subset \Omega$ be a complex submanifold. Then the restriction map
\[
 F \in \O(\Omega) \to F|_S \in \O(S)
\]
is a surjection.
\end{thm}
\begin{pf}
We take analytic polyhedra $\Delta_\nu \Subset \Omega$ and Oka maps
 $\phi_\nu: \bar\Delta_\nu (\Subset U_\nu) \to
 \overline{\PD}_\nu (\Subset \rP_\nu)$
 ($\nu=1,2, \ldots$) as in the proof of Theorem~\ref{thmchi}.
By Theorem~\ref{gsyzygy} (i) there is a finite generator system
$\{\underline{\sigma}_{\nu j}\}_{j=1}^{N_\nu}$ of
 $\sI\ang{S \cap \rP_\nu}$
on each $\overline\PD_\nu ~(\Subset \rP_\nu)$,
where $U_\nu$ is identified with $\phi_\nu(U_\nu)$.

Let $f \in \O(S)$ be any element.
 By Oka's J\^oku-Ik\^o Lemma~\ref{ojki} there
are $G_\nu \in \O(\overline{\PD}_\nu)$ 
with $G_\nu|_{\bar\Delta_\nu \cap S}=f|_{\bar\Delta_\nu \cap S}$
($\nu=1,2, \ldots$).

We set $F_1=G_1|_{\bar\Delta_1}$.
 Suppose that $F_\nu \in \O(\bar\Delta_\nu), 1 \leq \nu \leq \mu$,
are determined so that
\begin{align}
 F_\nu =f|_{\bar\Delta_\nu \cap S}, 
\quad
\| F_{\nu+1} - F_\nu\|_{\overline\PD_\nu} < \frac{1}{2^\nu},
\qquad 1 \leq \nu  \leq \mu -1.
\end{align}
For $\nu=\mu+1$ we first note that
$(G_{\mu+1}|_{\bar\Delta_\mu } - F_\mu)|_{\bar\Delta_\mu \cap S}=0$.
By Lemma~\ref{ojki} there is an element
$H_{\mu} \in \O(\overline\PD_{\mu})$ with
$H_{\mu}|_{\bar\Delta_\mu}=
G_{\mu+1}|_{\bar\Delta_\mu } - F_\mu$.
Since $\underline{H}_{\mu} \in \Gamma(\overline\PD_\nu,
 \sI\ang{S})$,
by Theorem~\ref{gsyzygy} (ii) there are $h_{\mu j} \in
 \O(\overline\PD_\mu)$, $1 \leq j \leq N_{\mu+1}$,
such that
\[
 H_{\mu} = \sum_{j=1}^{N_{\mu+1}} h_{\mu j} \cdot
 \sigma_{\mu+1 j} \quad \hbox{on } \overline\PD_\mu.
\]
Restricting this to $\bar\Delta_\nu$, we have
\[
 G_{\mu+1}|_{\bar\Delta_\mu}=F_\mu + \sum_{j=1}^{N_{\mu+1}} h_{\mu j} \cdot
 \sigma_{\mu+1 j}|_{\bar\Delta_\mu}.
\]
Approximating $h_{\mu j}$ sufficiently close by
 $\tilde{h}_{\mu j} \in
 \O(\Omega)$ on $\bar\Delta_\mu$ (Theorem~\ref{rwo}), and setting
\[
 F_{\mu+1} =  G_{\mu+1} - \sum_{j=1}^{N_{\mu+1}} \tilde{h}_{\mu j} \cdot
 \sigma_{\mu+1 j} \in \O(\bar\Delta_{\mu+1}),
\]
we have
\[
 F_{\mu+1}|_{\bar\Delta_{\mu+1} \cap S}=f|_{\bar\Delta_{\mu+1} \cap S},
\quad
\| F_{\mu+1} - F_\mu \|_{\bar\Delta_\mu} < \frac{1}{2^{\mu}}.
\]
Then the series
\[
 F=F_\mu + \sum_{\nu=\mu}^\infty (F_{\nu+1} -F_\nu)
\]
converges locally uniformly to the limit $F \in \O(\Omega)$ with $F|_S=f$.
\end{pf}

\begin{rmk}\rm
The above proof of Theorem~\ref{ext} works on Stein manifolds.
\end{rmk}

\subsection{Levi's (Hartogs' Inverse) Problem}\label{levip}
\subsubsection{Levi's (Hartogs' Inverse) Problem}
We first recall some basic terminologies (see, e.g., \cite{grr},
\cite{hor}, \cite{nog16} for more details).
Let $M$ be a connected complex manifold.

An upper semi-continuous function $\psi$ on $M$
 is said to be {\em plurisubharmonic} if every restriction
of $\psi$ to a 1-dimensional complex submanifold of
any holomorphic local chart of $M$ is subharmonic.
If $\psi$ is $C^2$-class, it is plurisubharmonic if and only if
the hermitian matrix
\[
\left(\frac{\del^2 \psi}{\del z_j \delbar z_k }\right)
\geq 0 \quad (\hbox{semi-posititive definite}),
\]
where $z=(z_j)$ is a holomorphic local coordinate system of $M$;
moreover, if
\[
\left(\frac{\del^2 \psi}{\del z_j \delbar z_k }\right)
\gg 0\quad (\hbox{posititive definite}),
\]
then $\psi$ is said to be {\em strongly plurisubharmonic}.

If $M$ carries a continuous plurisubharmonic exhaustion
 $\psi$\footnote{A functions $\psi$ is called an
exhaustion if $\{x \in M: \psi(x) < c\} \Subset M$
for all $c \in \R$.}, then $M$ is said to be {\em pseudoconvex}.
Note that there are several equivalent definitions of
being ``pseudoconvex'' (cf.\ \cite{guro}, \cite{grr}, \cite{hor}, \cite{nis},
\cite{nog16}, etc.).

A relatively compact domain $\Omega \Subset M$  is said to be
{\em strongly pseudoconvex} if for every boundary point
 $b \in \del \Omega$ there are a neighborhood $U$ of $b$ in $M$
and a strongly plurisubharmonic function $\phi$ on $U$
satisfying
\[
 \Omega \cap U=\{x \in U: \phi(x)<0\}.
\]
It is known that a strongly pseudoconvex domain is pseudoconvex.

If $M$ is connected and
 there is a locally finite holomorphic map $\pi: M \to \C^n$
with $\dim M=n$, then we call $M$ a Riemann domain,
in general; furthermore, if $\pi$ is locally (resp.\ non-) biholomorphic,
we call  $M$ a (resp.\ ramified) {\em unramified Riemann domain};
in this case, $M$ carries a Riemannian metric induced from the euclidean
one on $\C^n$ through $\pi$, so that $M$ satisfies the second
 countability axiom.

It had been known that a Stein manifold is pseudoconvex:
Levi's (Hartogs' Inverse) Problem had asked originally the converse for
univalent domains of $\C^n$.
K. Oka extended the problem for Riemann domains.
(There was  a necessity to do so (cf., e.g., \cite{nog16} \S5.1).)

\subsubsection{Oka's method}
Notice that Oka's J\^oku-Ik\^o Lemma~\ref{ojki}
 is sufficient to
deduce {\it Oka's Heftungslemma} which,
 together with a method of an integral equation and
the construction of
a plurisubharmonic exhaustion on a
 pseudoconvex unramified Riemann domain,  implies  Levi's (Hartogs' Inverse) Problem
 (cf.\ Oka \cite{okproc}, \cite{ok6},
 \cite{okp}, \cite{ok9},
Andreotti-Nara\-simhan \cite{an},   Nishino \cite{nis}):

\begin{thm}[Oka, 1941/42/43/53; cf.\ \S\ref{levisol}]\label{oka}
Let $M$ be a unramified Riemann domain.
If $M$ is pseudoconvex, then $M$ is Stein.
\end{thm}

\subsubsection{Grauert's method}
In 1958 H. Grauert \cite{gr58} gave another proof
 of Theorem \ref{oka} by
proving the finite dimensionality of the first cohomology of
coherent sheaves which was inspired
by the Cartan--Serre Theorem for coherent sheaves
 on compact analytic spaces.\footnote{
Cf.\ the footnote of \cite{gr58} p.~466.
The proof relies on L. Schwartz's finiteness theorem, whose rather
simple, short and complete proof is found in \cite{dem}
 and \cite{nog16} pp.~313--315.}
We shall observe that the Weak Coherence Theorem~\ref{weakc} suffices
for Grauert's method to prove Theorem~\ref{oka}.

We use the first \v{C}ech cohomology $H^1 ( \ast , \O)$.
The following immediately follows from Theorem \ref{thmchi}:
\begin{lem}\label{h1}
\begin{enumerate}
\item
If $\Omega$ be a holomorphically convex domain of $\C^n$,
then  $H^1(\Omega, \O)=0$.
\item
Let $M$ be a complex manifold with the second countability axiom, and let
 $\sU=\{U_\nu\}$ be a Stein covering of $M$ (i.e.,
every $U_\nu$ is open and Stein). Then we have
\[
 H^1(\sU, \O)\iso H^1(M, \O).
\]
\end{enumerate}
\end{lem}

Then we can apply Grauert's bumping method \cite{gr58} to prove:
\begin{thm}[Grauert]\label{gr}
Let $\Omega \Subset M$ be a relatively compact domain of a
complex manifold $M$ with strongly pseudoconvex boundary.
 Then we have
\[
  \dim_\C H^1(\Omega, \O) < \infty.
\]
\end{thm}

{\em Proof of Theorem \ref{oka}}:
(i)  The key of the proof is to show that a strongly pseudoconvex domain
$\Omega \Subset M$ is Stein. Let $b \in  \del \Omega$ be a boundary point.
By the definition of strong pseudoconvexity there are a neighborhood
$U$ of $b$ in $M$ and a quadratic polynomial $P_b(z_1, \ldots, z_n)$
such that $U \cap \bar\Omega \cap \{P_b=0\}=\{b\}$.
Then there is a little bit larger strongly pseudoconvex domain
 $\Omega_\epsilon \Supset\Omega$ such that $\Omega_\epsilon \cap \{P_b=0\}$
is closed in $\Omega_\epsilon$. We consider Cousin-I data
$(U_0, f_{\nu 0})$ and $(U_1, f_{\nu 1})$ of $\Omega_\omega$ such that
for $\nu=1,2, \ldots$,
\begin{equation}\nonumber
\begin{array}{ll}
 f_{\nu 0} =\dfrac{1}{P_b(z)^\nu}  &  \hbox{on } U_0 =
  \Omega_\epsilon\cap U, \\
 f_{\nu 1} =0   &  \hbox{on } U_1 = \Omega_\epsilon\setminus \{P_b=0\}.
\end{array}
\end{equation}
Theorem \ref{gr} applied to these Cousin-I data on
 $\Omega_\epsilon$ yields
 a meromorphic function $F_b$ on $\Omega_\epsilon$ such that
$F_b$ is holomorphic in $U_1$, and in $U_0$, $F_b$ is written as
\begin{equation}\label{bpole}
 F_b(z)=\frac{c_\nu}{P_b(z)^\nu}+ \cdots +\frac{c_1}{P_b(z)}
+\hbox{holomorphic term},
\end{equation}
where $c_\mu \in \C$ $(1 \leq \mu \leq \nu)$ and $c_\nu\not=0$.
Thus $F_b|_{\Omega} \in \O(\Omega)$ and $\Omega$ being
 holomorphically convex is deduced.

(ii) To show the holomorphic separation of $\Omega$, we take
two distinct points $Q_1, Q_2 \in \Omega$. We may assume that
$\pi(Q_1)=\pi(Q_2)=a \in \C^n$. Let $\phi(t), t \geq 0$, be
any affine linear curve with $\phi(0)=a$. Then there is
a unique lifting $\phi_j(t) \in \Omega, j=1,2$, such that
$\phi_j(0)=Q_j$ and $\pi \circ \phi_j(t)=\phi(t)$.
Since $\Omega$ is relatively compact, $\phi_j(t)$ hits the boundary
$\del \Omega$. We may assume that $\phi_1(t)$ hits $\del\Omega$
first with $t=T \in \R$, so that $\phi_j([0,T]) \subset \bar\Omega$
($j=1,2$) and $\phi_1(T) \in \del\Omega$. 
Note that $\phi_1(T) \not= \phi_2(T)$. With setting $b=\phi_1(T)$ we have by
\eqref{bpole} a meromorphic function $F_b$ in $\Omega_\epsilon$
which is holomorphic in $\Omega$.

We consider the Taylor expansions of $F_b$ at $Q_1$ and $Q_2$
in coordinates $(z_1, \ldots, z_n)$. Since $F_b$ has a pole at
$\phi_1(T)$ and no pole at $\phi_2(T)$, those two expansions must
be different. Therefore, there is some partial differential
operator
 $\del^\alpha=\del^{|\alpha|}/\del z_1^{\alpha_1} \cdots
\del z_n^{\alpha_n}$ with a multi-index $\alpha$
 such that
\[
 \del^\alpha F_b (Q_1) \not= \del^\alpha F_b(Q_2).
\]
Since $\del^\alpha F_b$ is holomorphic in $\Omega$,
this finishes the proof.
\qed

\begin{rmk}\rm\begin{enumerate}
\setlength{\itemsep}{-3pt}
\item
 The idea of the proof above is inspired
by Oka's unpublished paper in 1943 (\cite{oka43}).
Note that it is subtle how to deal with the holomorphic separation
of an unramified Riemann domain.
In H\"ormander \cite{hor} the holomorphic separation is included in the definition
 of Riemann domains. In \cite{nog16} Chap.~7 we used  Grauert's
Theorem \ref{gr} for a non-singular geometric ideal sheaves,
which can be also deduced from our Weak Coherence Theorem, but
then it involves a sheaf-cohomological argument.
In Gunning--Rossi \cite{guro} and in Nishino \cite{nis} they gave
their own proofs.
\item
The proof presented above provides a complete proof of Levi's
(Hartogs' Inverse) Problem without the sheaf cohomology theory nor
$L^2$-$\delbar$ method.
\end{enumerate}
\end{rmk}

\section{Historical remarks}\label{levisol}
Here, cf.\ \cite{nog16} Chap.~9 ``On Coherence'', and
cf., e.g.,  Lieb \cite{li} for the general background.

Oka's Theorem~\ref{oka} was first proved for univalent domains
 $\Omega \subset \C^2$
 by Oka \cite{okproc} (announcement)
in 1941, and  the full paper \cite{ok6} was published in 1942
with a comment of the validity for $n \geq 3$.

 In 1943 Oka proved Theorem~\ref{oka} for unramified Riemann domains of
general dimension $\geq 2$ in a series of research reports
of pp.~109 in total, sent to Teiji Takagi: The reports were
written in Japanese and unpublished (see \cite{okp}, \cite{oka43}).
  He remarked this fact three times, first
in his survey note \cite{oka49} (1949),
VIII \cite{ok8} (1951) and IX \cite{ok9} (1953).
The paper \cite{oka49} has not been referred very much, but
it should have a significant interest, for it was written
during the submission of VII \cite{ok7} before the publication;
 he had sent it to
H. Cartan 1948 one year ago (through the hands of S. Kakutani
and A. Weil),  and in return he had received
 Cartan's conjectural (or experimental)
 paper \cite{ca44} referred in \cite{oka49}.
 He surveyed the state of the
development of analytic function theory of several
variables at the time.

In the 1943 reports to Takagi he did not use Weierstrass'
Preparation Theorem, but he was writing a primitive form
of the notion of coherence and non-reduced structures
of analytic subsets; the study later led to the notion
``{\em id\'eaux de domaines ind\'etermin\'es}'', coherence in
 1948 (\cite{ok7}).
 The key of Oka's proof of Theorem~\ref{oka}
 was his ``{\it Heftungslemma}''.
In \cite{okproc} and  \cite{ok6} he proved Heftungslemma by Weil's integral,
but in 1943 (\cite{okp} no.~1, \cite{oka43})
 he replaced Weil's integral by simple Cauchy's integral,
 proving ``Oka's J\^oku-Ik\^o'' for unramified
Riemann domains.

In 1949  S. Hitotsumatsu \cite{hi} written in Japanese
 gave a proof of Oka's Heftungslemma by Weil's integral
to solve Levi's (Hartogs' Inverse) Problem in general dimension $n \geq 2$;
here he gave no argument of plurisubharmonic exhaustions
on pseudoconvex unramified Riemann domains, and so the result might hold
only for univalent  domains.

In 1953  Oka \cite{ok9} proved Theorem \ref{oka} above by making use of
his First and Second Coherence Theorems obtained in \cite{ok7}:
the Third Coherence Theorem was not used there.

In 1954  Bremermann \cite{br} and Norguet \cite{nor}
 independently proved Theorem \ref{oka}
for univalent domains $\Omega \subset \C^n$ with general $n \geq 2$,
 generalizing Oka's Heftungslemma
 by means of Weil's integral, similarly to Hitotsumatsu~\cite{hi}.

{\it Concluding Remark (Oka's Problem).} 
It is interesting to learn that Oka invented and proved three fundamental
coherence theorems by means of Weierstrass' Preparation Theorem
in order to treat the pseudoconvexity problem on singular 
ramified Riemann domains.
Levi's (Hartogs' Inverse) Problem for ramified Riemann domains has a counter-example
 (Forn{\ae}ss \cite{forn}),
 but in the same time there is a positive case for which
Levi's (Hartogs' Inverse) Problem is affirmative (\cite{nog17}).
Because of the above historical facts we may call
the following question

{\it Oka's Problem:  What is necessary and/or sufficient
for the validity of Levi's (Hartogs' Inverse) Problem on a ramified
Riemann domain $X$ (over $\C^n$)\,? :}

 This is open even when $X$ is non-singular.

\bigskip
\setlength{\baselineskip}{12pt}
\rightline{Graduate School of Mathematical Sciences}
\rightline{University of Tokyo (Emeritus)}
\rightline{Komaba, Meguro-ku, Tokyo 153-8914}
\rightline{Japan}
\rightline{e-mail: noguchi@ms.u-tokyo.ac.jp}

\bigskip
\rightline{(Revised 3/July/2018)}
\end{document}